\documentclass[12pt]{article}
\usepackage{longtable}
\usepackage{multirow}
\usepackage{scalefnt}
\usepackage{amsmath,amsthm,amsfonts,amssymb}
\usepackage{bm}
\usepackage{here}
\usepackage{caption}
\usepackage{setspace}
\usepackage[margin=15truemm]{geometry}
\usepackage{hyperref}
\date{Today}
\begin{document}
\begin{center}
{\LARGE\bf A family of Mordell curves with the rank of \\ at least three} \\
\vskip 1cm
\large
Seiji Tomita \\
\vskip 1cm
\end{center}

\begin{abstract}
We constructed a parametrized family of Mordell curves with the rank of at least three.
\end{abstract}

\vskip3\baselineskip

\section{Introduction}
\begin{equation}
y^2=x^3+d
\end{equation}

The elliptic curve represented by (1) is known as a Mordell curve.
Choudry and Zargar\cite{a} obtained a parametrized family of Mordell curves with the rank of at least three.
They obtained the solution using a quadratic system of equations.
Mina and Bacani\cite{b} constructed two infinite families with the rank of at least three.
They reduced the problem to an elliptic curve to get the third point.
In this paper, we present an elementary construction for a parametrized family of Mordell curves with the rank of at least three.
\begin{equation}
E: y^2=x^3+a^2
\end{equation}

\section{Main result}
We shall construct a parametrized family of Mordell curves with the rank of at least three.
Forcing a to be $x$-coodinates of the point on $E$. Then we get
$a^3+a^2=a^2(a+1)=v^2$ for some rational number $v$.
Then $a+1$ must be a square number. Let $a=m^2-1$ some rational number m.
Thus we get a first point
$P_1=(m^2-1,\ m(m^2-1))$ on the elliptic curve

$$y^2 = x^3+(m^2-1)^2.$$

Now, we try to increase the rank by forcing the x-coordinates of the second point $P_2$ to be $x=m-1$.

This holds when $m(m+3)=u^2$ for some rational number $u$.

A simple calculation shows that $m-1=-\cfrac{-5+4k}{(k-1)(k+1)}$.

Hence, the point $$P_2=\Bigl(-\frac{-5+4k}{(k-1)(k+1)},\ \frac{(-5+4k)(k-2)(2k-1)}{(k-1)^2(k+1)^2}\Bigr)$$

is a rational point on 

$$y^2 = x^3+\frac{(-5+4k)^2(2k^2-4k+3)^2}{(k-1)^4(k+1)^4}.$$

We try to increase the rank by forcing the x-coordinates of 3rd point $P_3$ to be $x=-\cfrac{2k^2-4k+3}{(k-1)(k+1)}$.

This holds when $-(2k^2+4k-7)$ is a perfect rational square number.

A simple calculation shows that $$x=\frac{(n^2-4n+22)(n^2+4n+6)}{4(4+n)(n-2)(n+1)}.$$

Hence, the point $$P_3=\Bigl(\frac{(n^2-4n+22)(n^2+4n+6)}{4(4+n)(n-2)(n+1)},\ \frac{(n^2+4n+6)(n^2+8n-2)(n^2-4n+22)(n^2+2n+10)}{16(4+n)^2(n-2)^2(n+1)^2}\Bigr)$$

is a rational point on 
\begin{equation}
y^2 = x^3+\cfrac{(n^2+4n+6)^2(34+n^2+8n)^2(2+n^2)^2(n^2-4n+22)^2}{256(4+n)^4(n-2)^4(n+1)^4},
\end{equation}

which contains the points,
\begin{align*}
P_1&=\Bigl(\frac{(n^2+4n+6)(34+n^2+8n)(2+n^2)(n^2-4n+22)}{16(4+n)^2(n-2)^2(n+1)^2},\\[6pt]
   & \quad \quad \frac{(n^2+2n+10)^2(34+n^2+8n)(n^2-4n+22)(2+n^2)(n^2+4n+6)}{64(4+n)^3(n-2)^3(n+1)^3}\Bigr),\\[6pt]
P_2&=\Bigl(-\frac{(34+n^2+8n)(2+n^2)}{4(4+n)(n-2)(n+1)},\ \frac{(34+n^2+8n)(n^2-4n-14)(2+n^2)(n^2+2n+10)}{16(4+n)^2(n-2)^2(n+1)^2}\Bigr),\\[6pt]
P_3&=\Bigl(\frac{(n^2-4n+22)(n^2+4n+6)}{4(4+n)(n-2)(n+1)},\ \frac{(n^2+4n+6)(n^2+8n-2)(n^2-4n+22)(n^2+2n+10)}{16(4+n)^2(n-2)^2(n+1)^2}\Bigr).
\end{align*}

Specialization to $n=3$ yields the elliptic curve

$$y^2 = x^3+\frac{142945242561}{157351936}.$$

The points $$P_1=(\frac{378081}{12544},\ \frac{236300625}{1404928}),\ P_2=(-\frac{737}{112},\ -\frac{313225}{12544}),\ P_3=(\frac{513}{112},\ \frac{397575}{12544})$$

have regulator $83.3621963770719$ by Sage\cite{c}, then the three points are independent,

 and Silverman's specialization theorem\cite{d} shows the rank of the family of Mordell curves
 
  given by the parametric solution of $(3)$ is at least three.

\section{High rank examples}

We calculated $(3)$ for $n<100$ and found several elliptic curves of the rank $4,5,6$, and $7$.

\begin{table}[hbtp]
\centering
  \caption{High rank examples}
  \label{table:data_type}
  \begin{tabular}{ccccc}
    \hline
   $rank$ & $n$   \\
    \hline \hline

$4$  & $1/3,$ & $1/7,$ & $1/8,$ & $1/9,$ $\cdots$ \\[5pt]
$5$  & $1/5,$ & $2/11,$ & $3/4,$ & $4/9,$ $\cdots$ \\[5pt]
$6$  & $1/14,$ & $2/3,$ & $5/22,$ & $7/3,$ $\cdots$ \\[5pt]
$7$  & $1/15,$ & $15/7,$ & $37/14,$ & $41/22,$ $\cdots$ \\
    \hline
  \end{tabular}
\end{table}


\begin{thebibliography}{99}

\bibitem{a} A. Choudhry and A. S. Zargar, A parametrised family of Mordell curves with a rational point of order 3,
            Notes on Number Theory and Discrete Mathematics, 26(1),2020. 

\bibitem{b}   Renz Jimwel S. Mina and Jerico B. Bacani,  Families of Mordell Curves with Non-trivial Torsion and Rank of at Least Three, 
International Conference on Mathematics and Computing, 2022

\bibitem{c} SAGE software, Available at \url{ http://sagemath.org.}

\bibitem{d} J. H. Silverman, Advanced topics in the arithmetic of elliptic curves, Springer, New York,
1994.

\end{thebibliography}
\end{document}